\documentclass{amsart}

\usepackage{fancyhdr}
\usepackage{lastpage}
\usepackage{stmaryrd,yhmath}
\usepackage{tikz}
\usepackage[all]{xy}

\newcommand{\bdis}{\begin{displaymath}}
\newcommand{\edis}{\end{displaymath}}
\newcommand{\be}{\begin{equation}}
\newcommand{\ee}{\end{equation}}
\newcommand{\mbb}{\mathbb}
\newcommand{\mcal}{\mathcal}

\theoremstyle{definition}

\theoremstyle{remark}

\newtheorem*{mydef11}{{\bf Theorem 1}}

\newtheorem*{mydef12}{{\bf Theorem 2}}

\newtheorem*{mydef2}{{\bf Definition}}

\newtheorem*{mydef51}{{\bf Lemma 1}}

\newtheorem*{mydef52}{{\bf Lemma 2}}

\newtheorem*{mydef53}{{\bf Lemma 3}}

\newtheorem*{mydef54}{{\bf Lemma 4}}

\newtheorem*{mydef55}{{\bf Lemma 5}}

\newtheorem*{mydef56}{{\bf Lemma 6}}

\newtheorem*{mydef57}{{\bf Lemma 7}} 

\newtheorem*{mydef58}{{\bf Lemma 8}}

\newtheorem*{mydefCHF1}{{\bf Complete Hybrid Formula 1}} 

\newtheorem*{mydefCHF2}{{\bf Complete Hybrid Formula 2}}

\numberwithin{equation}{section}

\begin{document}

\title{Jacob's ladders and new families of $\zeta$-kindred real continuous functions. }

\author{Jan Moser}

\address{Department of Mathematical Analysis and Numerical Mathematics, Comenius University, Mlynska Dolina M105, 842 48 Bratislava, SLOVAKIA}

\email{jan.mozer@fmph.uniba.sk}


\begin{abstract}
In this paper we obtain, by our method of crossbreeding in certain set of $\zeta$-factorization formulas, the corresponding complete hybrid formulas. These are playing the role of criterion for selection of new families of $\zeta$-kindred real continuous functions.  
\end{abstract}
\maketitle

\section{Introduction}

\subsection{}

In this paper we use the following notions: 
\begin{itemize} 
	\item[(A)] crossbreeding in the set of $\zeta$-factorization formulas; 
	\item[(B)] complete hybrid formula, 
	\item[(C)] definition of $\zeta$-kindred elements in the set of real continuous functions, 
\end{itemize}  
that we have introduced in the paper \cite{11}. 

Let us present a short survey of these notions. 

\begin{itemize}
	\item[(a)] We begin with the set of functions 
	\be \label{1.1} 
	\begin{split}
		& f_m(t)\in\tilde{C}_0[T,T+U],\ U=o\left(\frac{T}{\ln T}\right),\ T\to\infty \\ 
		& m=1,\dots,M,\ M\in\mbb{N}, 
	\end{split} 
	\ee 
	where $M$ is arbitrary and fixed. 
	
	\item[(b)] Next we obtain, by application of the operator $\hat{H}$ (introduced in the paper \cite{8}, (3.6)), the vector-valued functions 
	\be \label{1.2} 
	\begin{split}
		& \hat{H}f_m(t)=(\alpha_0^{m,k_m},\dots,\alpha_{k_m}^{m,k_m},\beta_1^{k_m},\dots,\beta_{k_m}^{k_m}), \\ 
		& m=1,\dots,M,,\ 1\leq k_m\leq k_0,\ k_0\in\mbb{N}, 
	\end{split} 
	\ee  
	where $k_0$ is arbitrary and fixed. Simultaneously, we obtain by our algorithm (for the short survey of this one see 
	\cite{8}, (3.1) -- (3.11)), also the following set of $\zeta$-factorization formulas 
	\be \label{1.3} 
	\begin{split} 
	& \prod_{r=1}^{k_m}
	\left|\frac{\zeta\left(\frac 12+i\alpha_r^{m,k_m}\right)}{\zeta\left(\frac 12+i\beta_r^{k_m}\right)}\right|^2\sim 
	E_m(U,T)F_m[f_m(\alpha_0^{m,k_m})], \\ 
	& m=1,\dots,M,\ L\to\infty. 
	\end{split} 
	\ee  
	
	\item[(c)] Further, we will suppose that we have obtained the following complete hybrid formula 
	\be \label{1.4} 
	\begin{split} 
	& \mcal{F}\left\{\prod_{r=1}^{k_1}(\cdots),\prod_{r=1}^{k_2}(\cdots),\dots,\prod_{r=1}^{k_M}(\cdots),
	F_1[f_1(\alpha_0^{1,k_1})],\dots,F_M[f_m(\alpha_0^{M,k_M})]\right\} \\ 
	& = 1+\mcal{O}\left(\frac{\ln\ln T}{\ln T}\right)\sim 1,\ T\to\infty 
	\end{split} 
	\ee  
	after the finite number of stages of crossbreeding (every memeber of (\ref{1.3}) is the participant on this process of 
	crossbreeding) in the set (\ref{1.3}) - that is: after the finite number of eliminations of the external functions 
	\be \label{1.5} 
	E_m(U,T),\ m=1,\dots,M 
	\ee  
	from the set (\ref{1.3}). 
	
	\item[(d)] Now we see that the complete hybrid formula (\ref{1.4}) expresses the functional dependence of the set of vector-valued functions (\ref{1.2}). Consequently, the back-projection of this functional dependence of the set (\ref{1.2}) into the generating set (\ref{1.1}) leaves us with the following (see \cite{1}) 
	
	\begin{mydef2}
		We will call the subset 
		\be \label{1.6} 
		\{f_1(t),\dots,f_M(t)\},\ t\in [T,T+U] 
		\ee  
		of the real continuous functions (comp. (\ref{1.1})), for which there is the complete hybrid formula (\ref{1.4}), as the family of $\zeta$-kindred functions. 
	\end{mydef2}
\end{itemize} 

\subsection{} 

In this paper we obtain the following new families of $\zeta$-kindred real continuous functions: 
\be \label{1.7} 
\begin{split} 
&  \left\{\frac{1}{\cos^2t},\frac{\sin^2t}{\cos^4t},\ (t-\pi L)^\Delta \right\},\\ 
& t\in [\pi L,\pi L+U],\ U\in (0,\pi/2-\epsilon],\Delta>0,\ L\to\infty, 
\end{split} 
\ee  
and  
\be \label{1.8} 
\begin{split} 
	&  \left\{\frac{1}{\cos^2t},\cos t,\ \cos^2t,\ \cos^3t,\ \sin^2t,\ (t-2\pi L)^\Delta \right\},\\ 
	& t\in [2\pi L,2\pi L+U],\ U\in (0,\pi/2-\epsilon],\Delta>0,\ L\to\infty.  
\end{split} 
\ee 

Let us remind we have introduced in the papers \cite{1} -- \cite{11} new notions in the theory of the Riemann zeta-function based on Jacob's ladders (see \cite{1}). The present paper contains new results in this direction. 

\section{The first class of lemmas} 

By making use of our algorithm for generating $\zeta$-factorization formulas (see \cite{8}, (3.1) -- (3.11)) we obtain the following results. 

\begin{mydef51}
For the function 
\be \label{2.1} 
f_1(t)=\frac{1}{\cos^2t}\in \tilde{C}_0[\pi L,\pi L+U],\ U\in (0,\pi/2-\epsilon] 
\ee  
there are the vector-valued functions 
\be \label{2.2} 
\begin{split} 
& (\alpha_0^{1,k_1},\dots,\alpha_{k_1}^{1,k_1},\beta_1^{k_1},\dots,\beta_{k_1}^{k_1}), \\ 
& 1\leq k_1\leq k_0,\ k_0\in \mbb{N} 
\end{split} 
\ee  
such that the following factorization formula 
\be \label{2.3} 
\prod_{r=1}^{k_1}
\left|\frac{\zeta\left(\frac 12+i\alpha_r^{1,k_1}\right)}{\zeta\left(\frac 12+i\beta_r^{k_1}\right)}\right|^2\sim 
\frac{\tan U}{U}\cos^2(\alpha_0^{1,k_1}),\ L\to\infty 
\ee 
holds true, where 
\be \label{2.4} 
\begin{split} 
& \alpha_r^{1,k_1}=\alpha_r(U,L,k_1;f_1),\ r=0,1,\dots,k_1, \\ 
& \beta_r^{k_1}=\beta_r(U,L,k_1),\ r=1,\dots,k_1, \\ 
& \pi L<\alpha_0^{1,k_1}<2\pi L+U \ \Rightarrow \ 0<\alpha_0^{1,k_1}-\pi L<U. 
\end{split} 
\ee 
\end{mydef51}  

\begin{mydef52} 
For the function 
\be \label{2.5} 
f_2(t)=\frac{\sin^2t}{\cos^4t}\in\tilde{C}_0[\pi L,\pi L+U],\  U\in (0,\pi/2-\epsilon] 
\ee  
\be \label{2.6} 
\begin{split} 
	& (\alpha_0^{2,k_2},\dots,\alpha_{k_2}^{2,k_2},\beta_1^{k_2},\dots,\beta_{k_2}^{k_2}), \\ 
	& 1\leq k_2\leq k_0,\ k_0\in \mbb{N} 
\end{split} 
\ee 
such that the following factorization formula 
\be \label{2.7} 
\prod_{r=1}^{k_2}
\left|\frac{\zeta\left(\frac 12+i\alpha_r^{2,k_2}\right)}{\zeta\left(\frac 12+i\beta_r^{k_2}\right)}\right|^2\sim 
\frac 13\frac{\tan^3 U}{U}\frac{\cos^4(\alpha_0^{2,k_2})}{\sin^2(\alpha_0^{2,k_2})},\ L\to\infty 
\ee 
holds true, where 
\be \label{2.8} 
\begin{split} 
	& \alpha_r^{2,k_2}=\alpha_r(U,L,k_2;f_2),\ r=0,1,\dots,k_2, \\ 
	& \beta_r^{k_2}=\beta_r(U,L,k_2),\ r=1,\dots,k_2, \\ 
	& 0<\alpha_0^{2,k_2}-\pi L<U. 
\end{split} 
\ee 
\end{mydef52}

\begin{mydef53} 
	For the function 
	\be \label{2.9} 
	\bar{f}_\Delta(t,L)=\bar{f}_\Delta(t)=(t-\pi L)^\Delta\in\tilde{C}_0[\pi L,\pi L+U],\  U\in (0,\pi/2),\ \Delta>0
	\ee   
	there are vector-valued functions 
	\be \label{2.10} 
	\begin{split} 
		& (\bar{\alpha}_0^{\Delta,\bar{k}_\Delta},\dots,\alpha_{\bar{k}_\Delta}^{\Delta,\bar{k}_\Delta},\beta_1^{\bar{k}_\Delta},\dots,\beta_{\bar{k}_\Delta}^{\bar{k}_\Delta}), \\ 
		& 1\leq \bar{k}_\Delta\leq k_0,\ k_0\in \mbb{N} 
	\end{split} 
	\ee 
	such that the following factorization formula 
	\be \label{2.11} 
	\prod_{r=1}^{\bar{k}_\Delta}
	\left|\frac{\zeta\left(\frac 12+i\bar{\alpha}_r^{\Delta,\bar{k}_\Delta}\right)}{\zeta\left(\frac 12+i\beta_r^{\bar{k}_\Delta}\right)}\right|^2\sim 
	\frac{1}{1+\Delta}
	\left(\frac{U}{\bar{\alpha}_0^{\Delta,\bar{k}_\Delta}-\pi L}\right)^\Delta,\ L\to\infty 
	\ee 
	holds true, where 
	\be \label{2.12} 
	\begin{split} 
		& \alpha_r^{\Delta,\bar{k}_\Delta}=\alpha_r(U,L,\bar{k}_\Delta;\bar{f}_\Delta),\ r=0,1,\dots,\bar{k}_\Delta, \\ 
		& \beta_r^{\bar{k}_\Delta}=\beta_r(U,L,\bar{k}_\Delta),\ r=1,\dots,\bar{k}_\Delta, \\ 
		& 0<\bar{\alpha}_0^{\Delta,\bar{k}_\Delta}-\pi L<U. 
	\end{split} 
	\ee 
(see \cite{9},\ (2.1) -- (2.5),\ $L\rightarrow\ \pi L$).
\end{mydef53}

\section{Theorem 1} 

\subsection{The first stage of the crossbreeding}  

It is the result of the crossbreeding between the $\zeta$-factorization formula (\ref{2.3}) and (\ref{2.7}):  
\be \label{3.1} 
\begin{split} 
& \prod_{r=1}^{k_1}
\left|\frac{\zeta\left(\frac 12+i\alpha_r^{1,k_1}\right)}{\zeta\left(\frac 12+i\beta_r^{k_1}\right)}\right|^2\sim \\ 
& \sim 
\frac{3}{U^2}\frac{\cos^6(\alpha_0^{1,k_1})\sin^2(\alpha_0^{2,k_2})}{\cos^4(\alpha_0^{2,k_2})}
\prod_{r=1}^{k_2}
\left|\frac{\zeta\left(\frac 12+i\alpha_r^{2,k_2}\right)}{\zeta\left(\frac 12+i\beta_r^{k_2}\right)}\right|^2,\ 
L\to\infty. 
\end{split} 
\ee  

\subsection{The second stage of the crossbreeding} 

Now, the crossbreeding between the formula (\ref{3.1}) and the formula (see (\ref{2.11}))
\be \label{3.2} 
U^\Delta\sim (1+\Delta)(\bar{\alpha}_0^{\Delta,\bar{k}_\Delta}-\pi L)^\Delta\prod_{r=1}^{\bar{k}_\Delta}
\left|\frac{\zeta\left(\frac 12+i\bar{\alpha}_r^{\Delta,\bar{k}_\Delta}\right)}{\zeta\left(\frac 12+i\beta_r^{\bar{k}_\Delta}\right)}\right|^2
\ee  
gives the following 
\begin{mydefCHF1}
\be \label{3.3} 
\begin{split} 
& \prod_{r=1}^{\bar{k}_\Delta}
\left|\frac{\zeta\left(\frac 12+i\bar{\alpha}_r^{\Delta,\bar{k}_\Delta}\right)}{\zeta\left(\frac 12+i\beta_r^{\bar{k}_\Delta}\right)}\right|^2\sim \\ 
& \sim \frac{3^{\Delta/2}}{1+\Delta}
\left[ 
\frac{|\cos^3(\alpha_0^{1,k_1})||\sin(\alpha_0^{2,k_2})|}
{(\bar{\alpha}_0^{\Delta,\bar{k}_\Delta}-\pi L)\cos^2(\alpha_0^{2,k_2})}
\right]^\Delta\times \\ 
& \times 
\left\{
\prod_{r=1}^{k_1}
\left|\frac{\zeta\left(\frac 12+i\alpha_r^{1,k_1}\right)}{\zeta\left(\frac 12+i\beta_r^{k_1}\right)}\right|^2
\right\}^{-3\Delta/2} \times 
\left\{
\prod_{r=1}^{k_2}
\left|\frac{\zeta\left(\frac 12+i\alpha_r^{2,k_2}\right)}{\zeta\left(\frac 12+i\beta_r^{k_2}\right)}\right|^2
\right\}^{\Delta/2}, \\ 
& 1\leq \bar{k}_\Delta,k_1,k_2\leq k_0,\ \Delta>0,\ L\to\infty. 
\end{split} 
\ee 
\end{mydefCHF1} 

\subsection{} 

Now, we obtain from (\ref{3.3}) by Definition the following 

\begin{mydef11} 
The subset 
\be \label{3.4} 
\begin{split} 
	&  \left\{\frac{1}{\cos^2t},\frac{\sin^2t}{\cos^4t},\ (t-\pi L)^\Delta \right\},\\ 
	& t\in [\pi L,\pi L+U],\ U\in (0,\pi/2-\epsilon],\Delta>0,\ L\to\infty 
\end{split} 
\ee  
is the family of $\zeta$-kindred elements in the class of real continuous functions. 
\end{mydef11}  

\section{The second class of lemmas} 

\subsection{} 

If we put in \cite{8}, (4.1) -- (4.10) 
\bdis 
L \rightarrow 2L,\ \mu=0 
\edis  
then we obtain the following result. 

\begin{mydef54} 
For the function 
\be \label{4.1} 
f_3(t)=\sin^2t\in\tilde{C}_0[2\pi L,2\pi L+U],\ U\in (0,\pi/2) 
\ee  
there are vector-valued functions 
\be \label{4.2} 
(\alpha_0^{3,k},\dots,\alpha_{k_3}^{3,k_3},\beta_1^{k_3},\dots,\beta_{k_3}^{k_3}),\ 1\leq k_3\leq k_0 
\ee  
such that the following factorization formula 
\be \label{4.3} 
\begin{split}
& \prod_{r=1}^{k_3} 
\left|\frac{\zeta\left(\frac 12+i\alpha_r^{3,k_3}\right)}{\zeta\left(\frac 12+i\beta_r^{k_3}\right)}\right|^2\sim \\ 
& \sim \left\{\frac 12-\frac 12\frac{\sin U}{U}\cos U\right\}
\frac{1}{\sin^2(\alpha_0^{3,k_3})},\ L\to\infty 
\end{split} 
\ee 
holds true, where 
\be \label{4.4} 
\begin{split} 
	& \alpha_r^{3,k_3}=\alpha_r(U,2L,k_3;f_3),\ r=0,1,\dots,k_3, \\ 
	& \beta_r^{k_3}=\beta_r(U,2L,k_3),\ r=1,\dots,k_3, \\ 
	& 0<\alpha_0^{3,k_3}-2\pi L<U. 
\end{split} 
\ee 
\end{mydef54}  

\begin{mydef55} 
For the function 
\be \label{4.5} 
f_4(t)=\cos^2t\in\tilde{C}_0[2\pi L,2\pi L+U],\ U\in (0,\pi/2) 
\ee  
there are vector-valued functions 
\be \label{4.6} 
(\alpha_0^{4,k_4},\dots,\alpha_{k_4}^{4,k_4},\beta_1^{k_4},\dots,\beta_{k_4}^{k_4}),\ 1\leq k_4\leq k_0 
\ee  
such that the following factorization formula 
\be \label{4.7} 
\begin{split}
	& \prod_{r=1}^{k_4} 
	\left|\frac{\zeta\left(\frac 12+i\alpha_r^{4,k_4}\right)}{\zeta\left(\frac 12+i\beta_r^{k_4}\right)}\right|^2\sim \\ 
	& \sim \left\{\frac 12+\frac 12\frac{\sin U}{U}\cos U\right\}
	\frac{1}{\cos^2(\alpha_0^{4,k_4})},\ L\to\infty 
\end{split} 
\ee 
holds true, where 
\be \label{4.8} 
\begin{split} 
	& \alpha_r^{4,k_4}=\alpha_r(U,2L,k_4;f_4),\ r=0,1,\dots,k_4, \\ 
	& \beta_r^{k_4}=\beta_r(U,2L,k_4),\ r=1,\dots,k_4, \\ 
	& 0<\alpha_0^{4,k_4}-2\pi L<U. 
\end{split} 
\ee 
\end{mydef55} 

\begin{mydef56} 
	For the function 
	\be \label{4.9} 
	f_5(t)=\frac{1}{\cos^2t}\in\tilde{C}_0[2\pi L,2\pi L+U],\ U\in (0,\pi/2-\epsilon] 
	\ee  
	there are vector-valued functions 
	\be \label{4.10} 
	(\alpha_0^{5,k_5},\dots,\alpha_{k_5}^{5,k_5},\beta_1^{k_5},\dots,\beta_{k_5}^{k_5}),\ 1\leq k_5\leq k_0 
	\ee  
	such that the following factorization formula 
	\be \label{4.11} 
	\begin{split}
		& \prod_{r=1}^{k_5} 
		\left|\frac{\zeta\left(\frac 12+i\alpha_r^{5,k_5}\right)}{\zeta\left(\frac 12+i\beta_r^{k_5}\right)}\right|^2\sim \\ 
		& \sim \frac{\sin U}{U}\frac{\cos^2(\alpha_0^{5,k_5})}{\cos U},\ L\to\infty 
	\end{split} 
	\ee 
	holds true, where 
	\be \label{4.12} 
	\begin{split} 
		& \alpha_r^{5,k_5}=\alpha_r(U,2L,k_5;f_5),\ r=0,1,\dots,k_5, \\ 
		& \beta_r^{k_5}=\beta_r(U,2L,k_5),\ r=1,\dots,k_5, \\ 
		& 0<\alpha_0^{5,k_5}-2\pi L<U, 
	\end{split} 
	\ee  
(comp. (\ref{2.1}) -- (\ref{2.4}) at $L\to\infty$). 
\end{mydef56} 

\subsection{} 

Next, we have the following results (see \cite{10}, (2.1) -- (2.6)). 

\begin{mydef57} 
	For the function 
	\be \label{4.13} 
	f_6(t)=\cos^3t\in\tilde{C}_0[2\pi L,2\pi L+U],\ U\in (0,\pi/2)
	\ee  
	there are vector-valued functions 
	\be \label{4.14} 
	(\alpha_0^{6,k_6},\dots,\alpha_{k_6}^{6,k_6},\beta_1^{k_6},\dots,\beta_{k_6}^{k_6}),\ 1\leq k_6\leq k_0 
	\ee  
	such that the following factorization formula 
	\be \label{4.15} 
	\begin{split}
		& \prod_{r=1}^{k_6} 
		\left|\frac{\zeta\left(\frac 12+i\alpha_r^{6,k_6}\right)}{\zeta\left(\frac 12+i\beta_r^{k_6}\right)}\right|^2\sim \\ 
		& \sim \left\{ \frac{\sin U}{U}-\frac{U^2}{3}\left(\frac{\sin U}{U}\right)^3 \right\}
		\frac{1}{\cos^3(\alpha_0^{6,k_6})},\ L\to\infty 
	\end{split} 
	\ee 
	holds true, where 
	\be \label{4.16} 
	\begin{split} 
		& \alpha_r^{6,k_6}=\alpha_r(U,L,k_6;f_6),\ r=0,1,\dots,k_6, \\ 
		& \beta_r^{k_6}=\beta_r(U,L,k_6),\ r=1,\dots,k_6, \\ 
		& 0<\alpha_0^{6,k_6}-\pi L<U. 
	\end{split} 
	\ee  
\end{mydef57}

\begin{mydef58} 
	For the function 
	\be \label{4.17} 
	f_7(t)=\cos t\in\tilde{C}_0[2\pi L,2\pi L+U],\ U\in (0,\pi/2)
	\ee  
	there are vector-valued functions 
	\be \label{4.18} 
	(\alpha_0^{7,k_7},\dots,\alpha_{k_7}^{7,k_7},\beta_1^{k_7},\dots,\beta_{k_7}^{k_7}),\ 1\leq k_7\leq k_0 
	\ee  
	such that the following factorization formula 
	\be \label{4.19} 
	\begin{split}
		& \prod_{r=1}^{k_7} 
		\left|\frac{\zeta\left(\frac 12+i\alpha_r^{7,k_7}\right)}{\zeta\left(\frac 12+i\beta_r^{k_7}\right)}\right|^2\sim \\ 
		& \sim  \frac{\sin U}{U}
		\frac{1}{\cos(\alpha_0^{7,k_7})},\ L\to\infty 
	\end{split} 
	\ee 
	holds true, where 
	\be \label{4.20} 
	\begin{split} 
		& \alpha_r^{7,k_7}=\alpha_r(U,L,k_7;f_7),\ r=0,1,\dots,k_7, \\ 
		& \beta_r^{k_7}=\beta_r(U,L,k_7),\ r=1,\dots,k_7, \\ 
		& 0<\alpha_0^{7,k_7}-2\pi L<U. 
	\end{split} 
	\ee  
\end{mydef58} 

\section{Theorem 2} 

\subsection{The first stage of the crossbreeding} 

The crossbreeding between the formula (\ref{4.3}) and (\ref{4.7}) gives the following formula 
\be \label{5.1} 
\begin{split}
& \frac{\cos^2(\alpha_0^{4,k_4})}{\cos U}\prod_{r=1}^{k_4} 
\left|\frac{\zeta\left(\frac 12+i\alpha_r^{4,k_4}\right)}{\zeta\left(\frac 12+i\beta_r^{k_4}\right)}\right|^2-
 \frac{\sin^2(\alpha_0^{3,k_3})}{\cos U}\prod_{r=1}^{k_3} 
\left|\frac{\zeta\left(\frac 12+i\alpha_r^{3,k_3}\right)}{\zeta\left(\frac 12+i\beta_r^{k_3}\right)}\right|^2 \\ 
& \sim 
\frac{\sin U}{U},\ L\to\infty. 
\end{split}
\ee 

\subsection{The second stage of the crossbreeding} 

The crossbreeding between the formula (\ref{4.15}) and (\ref{4.19}) gives the following formula 
\be \label{5.2} 
\begin{split}
	& \prod_{r=1}^{k_6} 
	\left|\frac{\zeta\left(\frac 12+i\alpha_r^{6,k_6}\right)}{\zeta\left(\frac 12+i\beta_r^{k_6}\right)}\right|^2 \sim \\ 
	& \sim 
	\frac{\cos(\alpha_0^{7,k_7})}{\cos^3(\alpha_0^{6,k_6})}
	\prod_{r=1}^{k_7} 
	\left|\frac{\zeta\left(\frac 12+i\alpha_r^{7,k_7}\right)}{\zeta\left(\frac 12+i\beta_r^{k_7}\right)}\right|^2- \\ 
	& - 
	\frac{U^2}{3\cos^3U\cos^3(\alpha_0^{6,k_6})}\times \\ 
	& 
	\left\{
	\cos^2(\alpha_0^{4,k_4})\prod_{r=1}^{k_4} 
	\left|\frac{\zeta\left(\frac 12+i\alpha_r^{4,k_4}\right)}{\zeta\left(\frac 12+i\beta_r^{k_4}\right)}\right|^2-
	\sin^2(\alpha_0^{3,k_3})\prod_{r=1}^{k_3} 
	\left|\frac{\zeta\left(\frac 12+i\alpha_r^{3,k_3}\right)}{\zeta\left(\frac 12+i\beta_r^{k_3}\right)}\right|^2
	\right\}^3, \\ 
	&  L\to\infty. 
\end{split}
\ee 

\subsection{The third stage of the crossbreeding} 

The crossbreeding between the formula (\ref{4.11}) and (\ref{5.1}) gives the following formula 
\be \label{5.3} 
\begin{split}
& \frac{1}{\cos^2U}\sim\frac{1}{\cos^2(\alpha_0^{5,k_5})}\prod_{r=1}^{k_5}
\left|\frac{\zeta\left(\frac 12+i\alpha_r^{5,k_5}\right)}{\zeta\left(\frac 12+i\beta_r^{k_5}\right)}\right|^2\times \\ 
& \left\{
\cos^2(\alpha_0^{4,k_4})\prod_{r=1}^{k_4}
\left|\frac{\zeta\left(\frac 12+i\alpha_r^{4,k_4}\right)}{\zeta\left(\frac 12+i\beta_r^{k_4}\right)}\right|^2- 
\sin^2(\alpha_0^{3,k_3})\prod_{r=1}^{k_3} 
\left|\frac{\zeta\left(\frac 12+i\alpha_r^{3,k_3}\right)}{\zeta\left(\frac 12+i\beta_r^{k_3}\right)}\right|^2
  \right\}^{-1}, \\ 
  & L\to \infty . 
\end{split}
\ee  

In the next stage of the crossbreeding we use the following variant of the formula (\ref{5.3}). 

\be \label{5.4} 
\begin{split}
& \frac{1}{\cos^3U}\sim \frac{1}{\cos^3(\alpha_0^{5,k_5})}
\left\{\prod_{r=1}^{k_5}\left|\frac{\zeta\left(\frac 12+i\alpha_r^{5,k_5}\right)}{\zeta\left(\frac 12+i\beta_r^{k_5}\right)}\right|^2\right\}^{3/2}\times \\ 
& 
\left\{
\cos^2(\alpha_0^{4,k_4})\prod_{r=1}^{k_4}
\left|\frac{\zeta\left(\frac 12+i\alpha_r^{4,k_4}\right)}{\zeta\left(\frac 12+i\beta_r^{k_4}\right)}\right|^2- 
\sin^2(\alpha_0^{3,k_3})\prod_{r=1}^{k_3} 
\left|\frac{\zeta\left(\frac 12+i\alpha_r^{3,k_3}\right)}{\zeta\left(\frac 12+i\beta_r^{k_3}\right)}\right|^2
\right\}^{-3/2}, \\ 
& L\to\infty . 
\end{split} 
\ee  

\subsection{The fourth stage of the crossbreeding} 

The crossbreeding between the formula (\ref{5.2}) and (\ref{5.4}) gives the following formula  
\be \label{5.5} 
\begin{split} 
& \prod_{r=1}^{k_6} 
\left|\frac{\zeta\left(\frac 12+i\alpha_r^{6,k_6}\right)}{\zeta\left(\frac 12+i\beta_r^{k_6}\right)}\right|^2\sim \\ 
& \sim \frac{\cos(\alpha_0^{7,k_7})}{\cos^3(\alpha_0^{6,k_6})}\prod_{r=1}^{k_7} 
\left|\frac{\zeta\left(\frac 12+i\alpha_r^{7,k_7}\right)}{\zeta\left(\frac 12+i\beta_r^{k_7}\right)}\right|^2- \\ 
& - \frac{U^2}{3}\frac{1}{\cos^3(\alpha_0^{5,k_5})\cos^3(\alpha_0^{6,k_6})}
\left\{
\prod_{r=1}^{k_5}\left|\frac{\zeta\left(\frac 12+i\alpha_r^{5,k_5}\right)}{\zeta\left(\frac 12+i\beta_r^{k_5}\right)}\right|^2
\right\}^{3/2}\times \\ 
& \left\{
\cos^2(\alpha_0^{4,k_4})\prod_{r=1}^{k_4}
\left|\frac{\zeta\left(\frac 12+i\alpha_r^{4,k_4}\right)}{\zeta\left(\frac 12+i\beta_r^{k_4}\right)}\right|^2- 
\sin^2(\alpha_0^{3,k_3})\prod_{r=1}^{k_3} 
\left|\frac{\zeta\left(\frac 12+i\alpha_r^{3,k_3}\right)}{\zeta\left(\frac 12+i\beta_r^{k_3}\right)}\right|^2
\right\}^{3/2}, \\ 
& L\to\infty . 
\end{split}
\ee  

\subsection{The fifth stage of the crossbreeding} 

In this stage we use the formula (\ref{2.11}) in the case 
\bdis 
L\rightarrow 2L;\ \pi L \rightarrow \pi 2L, 
\edis  
i.e. we use the following formula 
\be \label{5.6} 
\begin{split}
& U^2\sim (1+\Delta)^{2/\Delta}(\overset{*}{\bar{\alpha}}_0^{\Delta,\bar{k}_\Delta}-2\pi L)^2
\left\{
\prod_{r=1}^{\bar{k}_\Delta}
\left|\frac{\zeta\left(\frac 12+i\overset{*}{\bar{\alpha}}_r^{\Delta,\bar{k}_\Delta}\right)}{\zeta\left(\frac 12+i\overset{*}{\bar{\beta}}_r^{\bar{k}_\Delta}\right)}\right|^2, 
\right\}^{2/\Delta}, \\ 
& \Delta>0,\ L\to\infty , 
\end{split}  
\ee  
where (comp. (\ref{2.12})) 
\be \label{5.7} 
\begin{split} 
	& \overset{*}{\bar{\alpha}}_0^{\Delta,\bar{k}_\Delta}=\alpha_r(U,2L,\bar{k}_\Delta;\overset{*}{\bar{f}}_\Delta),\ r=0,1,\dots,\bar{k}_\Delta, \\ 
	& \overset{*}{\bar{\beta}}_r^{\bar{k}_\Delta}=\beta_r(U,2L,\bar{k}_\Delta),\ r=1,\dots,\bar{k}_\Delta, \\ 
	& \overset{*}{\bar{f}}_\Delta=(t-2\pi L)^\Delta\in \tilde{C}_0[2\pi L,2\pi L+U]. 
\end{split} 
\ee  
Finally, we obtain by the crossbreeding between the formulas (\ref{5.5}) and (\ref{5.6}) the 
following 

\begin{mydefCHF2}
\be \label{5.8} 
\begin{split}
& \prod_{r=1}^{k_6} 
\left|\frac{\zeta\left(\frac 12+i\alpha_r^{6,k_6}\right)}{\zeta\left(\frac 12+i\beta_r^{k_6}\right)}\right|^2\sim
\frac{\cos(\alpha_0^{7,k_7})}{\cos^3(\alpha_0^{6,k_6})}
\prod_{r=1}^{k_7} 
\left|\frac{\zeta\left(\frac 12+i\alpha_r^{7,k_7}\right)}{\zeta\left(\frac 12+i\beta_r^{k_7}\right)}\right|^2- \\ 
& - \frac{(1+\Delta)^{2/\Delta}}{3}
\frac{(\overset{*}{\bar{\alpha}}_0^{\Delta,\bar{k}_\Delta}-2\pi L)^2}{\cos^3(\alpha_0^{5,k_5})\cos^3(\alpha_0^{6,k_6})}
\left\{
\prod_{r=1}^{\bar{k}_\Delta}
\left|\frac{\zeta\left(\frac 12+i\overset{*}{\bar{\alpha}}_r^{\Delta,\bar{k}_\Delta}\right)}{\zeta\left(\frac 12+i\overset{*}{\bar{\beta}}_r^{\bar{k}_\Delta}\right)}\right|^2
\right\}^{2/\Delta} \times \\ 
& \times 
\left\{
\prod_{r=1}^{k_5}\left|\frac{\zeta\left(\frac 12+i\alpha_r^{5,k_5}\right)}{\zeta\left(\frac 12+i\beta_r^{k_5}\right)}\right|^2
\right\}^{3/2}\times \\ 
& 
\left\{
\cos^2(\alpha_0^{4,k_4})
\prod_{r=1}^{k_4}
\left|\frac{\zeta\left(\frac 12+i\alpha_r^{4,k_4}\right)}{\zeta\left(\frac 12+i\beta_r^{k_4}\right)}\right|^2-
\sin^2(\alpha_0^{3,k_3})
\prod_{r=1}^{k_3} 
\left|\frac{\zeta\left(\frac 12+i\alpha_r^{3,k_3}\right)}{\zeta\left(\frac 12+i\beta_r^{k_3}\right)}\right|^2
\right\}^{3/2}, \\ 
& 1\leq k_3,k_4,k_5,k_6,k_7,\bar{k}_\Delta\leq k_0,\ \Delta>0,\ L\to\infty.
\end{split}
\ee 
\end{mydefCHF2} 

\subsection{} 

Now, we obtain from (\ref{5.8}) by Definition the following 

\begin{mydef12} 
The set of functions 
\be \label{5.9} 
\begin{split}
& \left\{
\sin^2t, \cos^2t, \frac{1}{\cos^2t},\cos^3t,\cos t, (1-2\pi L)^\Delta
\right\}, \\ 
& t\in [2\pi L,2\pi L+U],\ U\in (0,\pi/2-\epsilon],\ \Delta>0,\ L\to\infty 
\end{split}
\ee 
is the family of $\zeta$-kindred elements in the class of real continuous functions. 
\end{mydef12}

\thanks{I would like to thank Michal Demetrian for their moral support of my study of Jacob's ladders.}

\end{document}